\documentclass[12pt]{article}
\usepackage{amssymb}
\usepackage{amsfonts}
\usepackage{amsmath}
\usepackage{amsmath}
\usepackage{verbatim}
\usepackage{color}
\usepackage[active]{srcltx}
\include{srctex}
\allowdisplaybreaks[4]
\newtheorem{theorem}{Theorem}[section]

\newtheorem{definition}[theorem]{Definition}

\newtheorem{problem}[theorem]{Problem}
\newtheorem{proposition}[theorem]{Proposition}
\newtheorem{remark}[theorem]{Remark}

\let\Section=\section
\def\section{\setcounter{equation}{0}\Section}
\newenvironment{proof}[1][Proof]{\textbf{#1.} }{\ \rule{0.5em}{0.5em}}

\def\RR{\mathbb{R}}
\def\NN{\mathbb{N}}
\def\EE{\mathbb{E}}
\def\SS{\mathbb{S}}
\def\PP{\mathbb{P}}

\def\cF{{\cal F}}

\def\si{{\sigma}}

\def\si{{\sigma}}

\def \eref#1{\hbox{(\ref{#1})}}

\begin{document}

\title{Gaussian estimates of the density for systems of non-linear stochastic heat equations
\thanks{ Supported by the NNSF of China (No.: 11271169) and the Priority Academic Program Development of Jiangsu Higher Education Institutions.\newline
Keywords and phrases: Stochastic heat equations, Malliavin calculus,
spatially homogeneous covariances, Gaussian density. }}

\author{{Yinghui Shi}$$\footnote{E-mail: shiyinghui@jsnu.edu.cn}\;
\  \ \ {Xiaobin Sun}$$\footnote{E-mail: xbsun@jsnu.edu.cn}\;
\\
 \small   School of Mathematics and Statistics, Jiangsu Normal University, Xuzhou 221116, China.}\,

\date{\today}
\maketitle

\begin{abstract}
In this paper we consider a system of non-linear stochastic heat equations on $\mathbb{R}^d$ driven by a Gaussian noise which is
white in time and has a homogeneous spatial covariance. Under some suitable regularity and non degeneracy conditions, the smoothness of the joint density of the solution for this system has been studied by E. Nualart in \cite{NuaE}. The purpose of this paper is further to study the lower and upper bounds of the density. The main tools are the Malliavin calculus and
the method developed by Kohatasu-Higa in \cite{K} or E. Nualart and Quer-Sardanyons in \cite{NuaEQ}.

\end{abstract}

\section{Introduction}
Consider the stochastic partial differential equations
\begin{equation}\label{Eq}
\frac{\partial u_i}{\partial t}(t,x)=\frac{1}{2}\frac{\partial^2 u_i }{\partial x^2}(t,x)+b_i(u(t,x))+\sum^{q}_{j=1}\sigma_{ij}(u(t,x))\dot{W}^j(t,x),\quad i=1,2\ldots,m,
\end{equation}
with vanishing initial conditions, $x\in \mathbb{R}^d$, $u=(u_1,\ldots, u_m)$.
Here $\sigma_{ij}, b_i: \mathbb{R}^m\rightarrow \mathbb{R}$ are globally Lipschitz functions, which are the entries of a $m\times q$ matrix $\sigma$ and a $m$-dimensional vector $b$, $\sigma=(\sigma_1,\ldots, \sigma_m)$ and $b=(b_1,\ldots, b_m)$. The perturbation $\dot{W}(t,x)=(\dot{W}^1(t,x),\ldots, \dot{W}^q(t,x))$ is a $q$-dimensional Gaussian noise which is white in time and with a spatially homogeneous covariance $f$, that is,
$$
\mathbb{E}[\dot{W}^i(t,x)\dot{W}^j(s,y)]=\delta(t-s)f(x-y)\delta_{ij},
$$
$\delta(\cdot)$ denotes the Dirac delta function, $\delta_{ij}$ is the Kronecker symbol, and $f$ is a positive continuous function on $\mathbb{R}^d\backslash \{0\}$.


A mild solution of Eq.\eref{Eq} can be formulated by using the Green kernel $\Gamma(t,x)$ associated with the operator $L=\frac{\partial}{\partial t}-\frac{1}{2}\frac{\partial^2}{\partial x^2}$ (see Definition  \ref{def mild sol}).
This requires the notion of stochastic integral introduced by Walsh in  \cite{Wal}.

The Malliavin calculus (see \cite{NuaD} or \cite{San}) is a powerful tool to study the existence and smoothness of the density for the solutions of SPDEs (see \cite{BP, HHNS, MNQ, MS, NuaE, NQ, QS1, QS2}), also the lower and upper bounds for the density (see \cite{GMEN, K, NuaEQ}). For the proof of the upper bound, a well-known technique is to get the expression of the density by integration by parts formula, then H\"older's inequality and exponential martingale inequality imply the estimate of the Malliavin norms for the derivative and the Malliavin matrix (see for instance \cite{GMEN} or \cite{NuaEQ}). However, the proof of the lower bound is more difficult. For a uniformly hypoelliptic diffusion with smooth drift, Kusuoka and Stroock obtained a Gaussian type lower bound in \cite{KS}. Kohatsu-Higa in \cite{K} extended the  results of Kusuoka and Stroock to general random variables on Wiener space, and studied one-dimensional stochastic  heat equation on $[0,1]$ driven by space-time white noise. Later on, E. Nualart and Quer-Sardanyons established the lower and upper bounds for the solution of the stochastic heat equation driven by a Gaussian noise with white in time and spatially homogeneous covariance in \cite{NuaEQ}.

In \cite{NuaE}, E. Nualart has studied the smoothness of the joint density of systems of non-linear spatially homogeneous SPDEs, which include Eq.\eref{Eq},
by using Malliavin calculus techniques. The aim of this paper is a further study of the lower and upper bounds of the density.
This result extends previous work of E. Nualart and Quer-Sardanyons in \cite{NuaEQ}.

The paper is organized as follows. Some preliminaries are given in Section 2. Section 3 is devoted to prove the lower and upper bounds for the density for solution $u(t,x)=(u_1(t,x),\ldots,u_m(t,x))$, following the general criterion established in \cite{K} or \cite{NuaEQ} .  Finally, the results are applied to the spatial covariances given by the Riesz, Bessel and fractional kernels.

\section{Preliminaries}\label{preliminaries}

Consider a  non-negative and non-negative definite function $f$ which is continuous on $\RR^d \setminus \{0\}$.
We  assume that $f$ is the Fourier transform
of a non-negative tempered measure $\mu$ on $\RR^d$ (called the
spectral measure of $f$). That is, for all $\varphi$ belonging to
the space $\mathcal{S}(\RR^d)$ of rapidly decreasing $\mathcal{C}^{\infty}$
functions on $\RR^d$
\begin{equation}\label{def of spectral measure mu}
\int_{\RR^d}f(x)\varphi(x)dx=\int_{\RR^d}\mathcal{F}\varphi(\xi)\mu(d\xi),
\end{equation}
and assume the following condition:
\begin{equation}
\Phi(T):=\int^{T}_{0}\int_{\RR^d}|\cF\Gamma(r)(\xi)|^2\mu(d\xi)dr<\infty\,,\label{2.2}
\end{equation}
which is used to prove the unique solution for Equation \eref{Eq}. Here we denote by $\mathcal{F}\varphi$ as the Fourier transform of
$\varphi\in \mathcal{S}(\RR^d)$, denoted by
$\mathcal{F}\varphi(\xi)=\int_{\RR^d}\varphi(x)e^{- i\xi\cdot x}dx$. Moreover, it has been proved that condition \eref{2.2} is equivalent to
\begin{equation}
\int_{\RR^d}\frac{1}{1+|\xi|^2}\mu(d\xi)<\infty\,.\label{2.3}
\end{equation}
We need a slightly stronger condition than \eref{2.3} in order to prove our main result.

\vspace{3mm}
$({\bf H_{\eta}})$ For some $\eta\in(0,1)$, it holds:
$$
\int_{\RR^d}\frac{1}{(1+|\xi|^2)^{\eta}}\mu(d\xi)<\infty\,.
$$

Suppose that $(\Omega, \mathcal{F}, \mathbb{P})$ is a complete probability space. For $T>0$, let $\mathcal{C}_0^{\infty}([0,T]\times \RR^{d})$ be the space of smooth
 functions with compact support on $[0,T] \times \RR^d$.  Consider a  family of zero mean Gaussian  random variables
$W=\{W^j(\varphi), j=1,\ldots, q, \varphi \in \mathcal{C}_0^{\infty}([0,T]\times\RR^{d})\}$ with covariance
\begin{equation}\label{cov}
\EE(W^i(\varphi)W^j(\psi))=\delta_{ij}\int_0^{T}\int_{\RR^d}\int_{\RR^d}\varphi(t,x)f(x-y)\psi(t,y)dxdydt\,.
\end{equation}
Using Fourier transform, \eref{cov} can also be written as
\begin{equation}
\EE (W^i(\varphi)W^j(\psi))=\delta_{ij}\int_0^T\int_{\RR^d} \mathcal{F}\varphi(t)(\xi) \overline {\mathcal{F}\psi(t)(\xi)}\mu(d\xi)dt\,,\nonumber
\end{equation}
where $\overline {\mathcal{F}\psi}$ is the complex conjugate of $\mathcal{F}\psi$.



Let $\mathcal{H}^q$ be the Hilbert space which is the completion of $\mathcal{C}_0^{\infty}(\mathbb{R}^d; \RR^q)$ with the inner product
\begin{equation} \label{def H}
\langle\varphi,\psi\rangle_{\mathcal{H}^q}=\sum^{q}_{l=1}\int_{\RR^d}dx\int_{\RR^d}dy\varphi_l(x)f(x-y)\psi_l(y)=\sum^{q}_{l=1}\int_{\RR^d} \mathcal{F}(\varphi_l)(\xi) \overline {\mathcal{F}(\psi_l)(\xi)}\mu(d\xi)
\end{equation}
for $\varphi,\psi \in \mathcal{C}_0^{\infty}(\RR^d; \RR^q)$. Notice that $\mathcal{H}^q$ may contain distributions.


The Gaussian family $W$ can be extended to the space $\mathcal{H}^q_T=L^2([0,T];\mathcal{H}^q)$ by $W(g)=\sum^q_{i=1}W^{i}(g_i)$, $g\in\mathcal{H}^q_T$.  It is obvious that ${\bf 1}_{[0,t]}h$ belongs to
$\mathcal{H}^q_T$. Put $W_t(h)=W({\bf 1}_{[0,t]}h)$ for $t\ge 0$ and $h \in \mathcal{H}^q$, we have that $W=\{ W_t, 0\leq t\leq T\}$ is a cylindrical Wiener
process in the Hilbert space $\mathcal{H}^q$ (see \cite{DJ}). That is, for any $h \in \mathcal{H}^q$,
$\{W_t(h), 0\leq t\leq T\}$ is a Brownian motion with variance
$t\| h\| ^2_{\mathcal{H}^q}$, and
\begin{equation*}
\EE(W_t(h)W_s(g))=(s\wedge t)\langle h,g\rangle_{\mathcal{H}^q} .
\end{equation*}

Let $(\mathcal{F}_t)_{t\geq 0}$ be the $\sigma$-filtration generated by the random
variables $\{W_s(h), h\in \mathcal{H}^q, 0\leq s \leq t\}$ and the $\mathbb{P}$-null
sets.
Define the stochastic integral for an $\mathcal{H}^q$-valued $\mathcal{F}_t$-predictable
process $g \in L^2(\Omega\times[0,T];\mathcal{H}^q)$ with respect to the
cylindrical Wiener process $W$ as
\begin{equation*}
\int^T_0 g(t)dW_t:=\int_0^T\int_{\RR^d}g(t,x)W(dt,dx):=\sum^q_{j=1}\int_0^T\int_{\RR^d}g_j(t,x) W^j(dt,dx),
\end{equation*}
then we have the isometry property
\begin{equation}\label{isometry property}
\EE\left|\int^T_0 g(t) dW_t\right|^2=\EE\int_0^T \| g(t)\|^2_{\mathcal{H}^q} dt.
\end{equation}

Using the notion of the above stochastic integral, we introduce the
following definition:

\begin{definition}\label{def mild sol}
A $\RR^m$-valued adapted stochastic process $\{u(t,x)=(u_1(t,x),\ldots, u_m(t,x)), (t,x)\in [0,T]\times \RR^d\}$ is a mild solution of Eq.\eref{Eq} if for all $0\leq t\leq T$, $x\in \RR^d$, $i=1,\ldots,m$,
\begin{eqnarray*}
 u_i(t, x)
&=&\sum^{q}_{j=1}\int_0^t\int_{\RR^d}\Gamma(t-s,x-y)\si_{ij}(u(s,y)) W^j(ds,dy)\nonumber\\
&&+\int_0^t\int_{\RR^d}b_i(u(s,y))\Gamma(t-s,x-y)dyds, \,\ \ \ \
\mathbb{P}-a.s.,
\end{eqnarray*}
where $\Gamma(t,x)=(2\pi t)^{-\frac{d}{2}}\exp{\{-\frac{|x|^2}{2t}\}}$ is the fundamental solution to $\frac{\partial u}{\partial t}-\frac{1}{2}\frac{\partial^2 u}{\partial x^2}=0$.
\end{definition}

Now, we state the existence, uniqueness and H\"older continuity of the solution for Eq.\eref{Eq}, which have been showed in \cite[Sections 2 and 3]{NuaE}.

\begin{theorem}\label{existence and uniqueness THM}
Assume condition \eref{2.2} is satisfied, then there exists a unique mild solution $u$ to Eq.\eref{Eq} such that for all $p\geq 1$ and $T>0$,
\begin{equation}
\sup_{(t,x) \in [0,T]\times\mathbb{R}^d}\EE|u(t,x)|^p< +\infty.
\label{moments estimate}
\end{equation}
Furthermore, if condition $({\bf H_{\eta}})$ holds, then for all $\gamma_1\in(0, \frac{1-\eta}{2})$, $s, t\in [0, T]$, $x\in\RR^d$ and $p>1$,
\begin{eqnarray}
\EE |u(t,x)-u(s,x)|\leq C_{p, T}|t-s|^{\gamma_1 p}\label{gamma_1}
\end{eqnarray}
and for all $\gamma_2\in(0, 1-\eta)$, $t\in [0, T]$, $x,y\in\RR^d$ and $p>1$,
\begin{eqnarray}
\EE |u(t,x)-u(t,y)|\leq C_{p, T}|x-y|^{\gamma_2 p}\label{gamma_2}
\end{eqnarray}
for some constant $C_{p, T}>0$.
\end{theorem}

Next we recall some concepts of Malliavin calculus which is used to prove the main results. Notice that $\{W(h), h \in \mathcal{H}^q_T\}$
is a centered Gaussian process and $\EE (W(h_1)W(h_2))=\langle h_1, h_2 \rangle_{\mathcal{H}^q_T}, h_1, h_2 \in \mathcal{H}^q_T$,
then we can develop a Malliavin calculus (see \cite{NuaD}). The Malliavin derivative is denoted by $D$, which is a closed operator on $L^2(\Omega)$
and takes the value in $L^2(\Omega ; \mathcal{H}_T^{q})$. For any integer $k\geq 1$ and $p\ge 2$, denote the domain of the iterated derivative
$D^k$ by $\mathbb{D}^{k, p}$ and $\mathbb{D}^{\infty}:=\cap_{p \geq 1} \cap_{k \geq 1} \mathbb{D}^{k,p}\,.$ The space $\mathbb{D}^{k,p}$ also is the completion of the set of smooth functionals with respect to seminorm
  \[
\|F\|_{k, p}=\left\{\EE[|F|^p]+\sum^{k}_{j=1}\EE[\|D^{j}F\|_{(\mathcal{H}^{q}_T)^{\otimes j}}]\right\}^{\frac{1}{p}}.
\]
For any $X\in \mathbb{D}^{1, 2}$ and some fixed $r  \ge 0$, $DX(r,*)$ is an element of $\mathcal{H}^q$, which will be denoted by $D_{r,*}X$.

We define the Malliavin matrix of a $m$-dimensional random vector $X=(X_1,\ldots, X_m)\in (\mathbb{D}^{1,2})^m$ by $M_X=(\langle DX_i, DX_j\rangle_{\mathcal{H}_T})_{1\leq i, j\leq k}$.
We will say that a random vector $X$ whose components are in $\mathbb{D}^{\infty}$ is non-degenerate if $(\det M_X)^{-1}\in \cap_{p\geq 1}L^p(\Omega)$. It is well-known that a non-degenerate random vector has a smooth density (see \cite[Proposition 2.1.5]{NuaD}).

As did in \cite{K}, we set $\mathcal{H}^q_{s, t}= L^{2}([s, t]; \mathcal{H}^q)$ and $\|\cdot\|_{s, t}:=\|\cdot\|_{\mathcal{H}^q_{s, t}}$. For any integer $k\geq 1$ and $p> 1$, we define the seminorm:
\[
\|F\|^{s, t}_{k, p}=\left\{\EE_s[|F|^p]+\sum^{k}_{j=1}\EE_s[\|D^{j}F\|_{(\mathcal{H}^{q}_{s, t})^{\otimes j}}]\right\}^{\frac{1}{p}},
\]
where $\EE_s[\cdot]=\EE[\cdot| \cF_s]$. We also write $\mathbb{P}_s\{\cdot\}=\mathbb{P}\{\cdot|\mathcal{F}_s\}$. Completing the space of smooth functionals with respect to this seminorm, we obtain the space $\mathbb{D}^{k,p}_{s,t}$. We say that $F\in \overline{\mathbb{D}}^{k,p}_{s,t}$  if $F\in \mathbb{D}^{k,p}_{s,t}$ and $\|F\|^{s, t}_{k, p}\in \cap_{q\geq 1}L^{q}(\Omega)$, and we set $\overline{\mathbb{D}}^{\infty}_{s,t}:=\cap_{k\geq 1}\cap_{p\geq 1}\overline{\mathbb{D}}^{k,p}_{s,t}$. Furthermore, we define the conditional Malliavin covariance matrix associated to an $m$-dimensional random vector $X=(X^1,\ldots, X^m)$ by $M^{s,t}_{X}:=(\langle DX^i, DX^j\rangle_{\mathcal{H}^q_{s, t}})_{1\leq i,j\leq m}$.

The next result is the $q$-dimensional extension of \cite[Proposition 6.1]{NQ} and \cite[Lemma 3.4]{NuaEQ}. The proof is omitted because it follows exactly the same arguments.

\begin{proposition}
 Assume that the coefficient $\sigma$, $b$ are smooth functions with bounded partial derivatives of order greater than or equal to one. Then, for any $(t,x)\in [0, T]\times \mathbb{R}^d$, the random variable $u_i(t,x)$ belongs to the space $\mathbb{D}^{\infty}$, for all $i=1,\ldots, m$. Moreover, the derivative $Du_i(t,x)$ is an $\mathcal{H}^{q}_T$-valued process which satisfies the following linear stochastic differential equation:
\begin{eqnarray}
D_{r, \ast}u_i(t, x)&=&\Gamma(t-r, x-\ast)\sigma_{i}\left(u(r,\ast)\right)\nonumber\\
&&+\sum^{q}_{l=1}\int_r^t\int_{\RR^d}\Gamma(t-s,x-y)D_{r, \ast}(\sigma_{il}(u(s, y)))W^{l}(ds,dy)\nonumber\\
&&+\int_0^t\int_{\RR^d}D_{r, \ast}(b_i(u(s,y)))\Gamma(t-s,x-y)dyds
\end{eqnarray}
for all $r\in [0,t]$, and $D_{r, \ast}u_i(t, x)=0$, for all $r>t$.

Furthermore, for $0\leq a<b\leq T$ and $p\geq 1$, there exits a positive constant $C=C(a, b)$ such that for all $\delta\in(0, b-a]$:
\begin{eqnarray}
\sup_{(t,x)\in [b-\delta, b]\times\RR^d}\EE_a\|D^m u_i(t,x)\|^{2p}_{(\mathcal{H}^q_{b-\delta, b})^{\otimes m}}\leq C(\Phi(\delta))^{mp}, \quad a.s.,\label{2.7}
\end{eqnarray}
where $\Phi(\delta)$ is the one in \eref{2.2}.
\end{proposition}

In order to prove the existence of the smooth density of $u(t,x)$, we need the following conditions:

\medskip
({\bf H1}) There exists $\beta >0$  such that for all $\varepsilon\in (0,1]$,
\begin{equation}
C \varepsilon^{\beta} \leq \int^{\varepsilon}_{0}\int_{\RR^d}|\cF\Gamma(r)(\xi)|^2\mu(d\xi)dr\nonumber
\end{equation}
for some constant $C> 0$.

\medskip
({\bf H2}) Let $\beta$ be given in hypothesis ({\bf H1}) and  $\gamma_1$ and $\gamma_2$ be given in \eref{gamma_1} and \eref{gamma_2}.

(i)
The function $\Psi(t,x):=|x|^{\gamma_2} \Gamma(t,x)$ satisfies $\int_0^T \int_{\RR^d} | \mathcal{F} \Psi(t) (\xi) |^2 \mu(d\xi) dt<\infty$  and
there exists $\beta_1>\gamma_2 \vee \beta$ such that for all $\varepsilon\in [0,1]$  satisfying
\begin{equation}\label{hyp 4-3}
\int_0^{\varepsilon} \langle  \Psi(r,*), \Gamma(r,*)\rangle_{\mathcal{H}}dr\leq C\varepsilon^{\beta_1},
\end{equation}
for some positive constant $C$.

(ii)  There exists $\beta_2  > \gamma_1 \vee \beta$ such that, for all $\varepsilon\in [0,1]$,
\begin{equation}\label{hyp 4-1}
\int_0^{\varepsilon} r^{\gamma_1}\int_{\RR^d}|\cF\Gamma(r)(\xi)|^2\mu(d\xi)dr \leq C\varepsilon^{\beta_2}
\end{equation}
for some positive constant $C$.

We have the following theorem (see \cite[Theorem 4.1]{NuaE}):
\begin{theorem} \label{THM smooth density}
Assume conditions $({\bf H_{\eta}})$, $({\bf H1})$ and $({\bf H2})$ hold, and the coefficients $\sigma$, $b$ are smooth functions with bounded partial derivatives of order greater than or equal to one. Then for all $(t,x)\in(0, T]\times \RR^d$, the law of the random vector $u(t,x)$ admits a smooth density $p_{t,x}(\cdot)$ on $\Sigma:=\{y\in\RR^m: \sigma_1(y),\ldots,\sigma_q(y) \ \ \text{span} \ \ \RR^m\}$.
\end{theorem}

\section{Lower and upper bound for the density}

In this section, we shall study the lower and upper bounded of the density $p_{t,x}(\cdot)$.
As the argument in \cite[Lemma 3.1]{MMS}, the condition \eref{2.3} implies that there exists positive constant $C_1$ such that
\begin{eqnarray}
C_1(t-s)\leq \int^{t}_{s}\int_{\RR^d}|\cF\Gamma(r)(\xi)|^2\mu(d\xi)dr,~~~~0\leq s< t\leq T.\label{ie.4.1}
\end{eqnarray}
Furthermore, if condition $({\bf H_{\eta}})$ holds, there exists positive constant $C_2$ such that
\begin{eqnarray}
\int^{t}_{0}\int_{\RR^d}|\cF\Gamma(r)(\xi)|^2\mu(d\xi)dr\leq C_2 t^{1-\eta}.\label{ie.4.2}
\end{eqnarray}

\begin{remark}
Similar to the argument in \cite[Remark 3.1]{NuaEQ}, the estimate \eref{ie.4.1} will play an important role in the proof of the lower bound. This has prevented us from considering the other type of SPDEs, such as the stochastic wave equation. Actually, we do not have a kind of time homogeneous lower bound of the form \eref{ie.4.1} for stochastic wave equation.
\end{remark}

In order to obtain the lower bound, we need more conditions on the coefficients $\sigma$ and $b$:

\medskip
({\bf H3}) Assume that $b_i$ are bounded, for any $i=1,\ldots, m$£¬ and there exist positive constants $C_1$ and $C_2$, such that for all $\xi\in\RR^m$,
\begin{eqnarray}
C_1|\xi|^2\leq\inf_{x,y\in\RR^d}\sum^m_{i,j=1}\sum^q_{k=1}\sigma_{ik}(x)\sigma_{jk}(y)\xi_i\xi_j \label{lower of sigma}
\end{eqnarray}
and
\begin{eqnarray}
\sup_{x,y\in\RR^d}\sum^m_{i,j=1}\sum^q_{k=1}\sigma_{ik}(x)\sigma_{jk}(y)\xi_i\xi_j \leq C_2|\xi|^2.\label{upper of sigma}
\end{eqnarray}

\begin{remark} If $m=1$ and $q=1$, then \eref{lower of sigma} and \eref{upper of sigma} is equivalent to $C_1\leq |\sigma(x)|\leq C_2$, for any $x\in\mathbb{R}$, which is the condition on $\sigma$ in \cite[Theorem 1.1]{NuaEQ}.
\end{remark}

The main theorem of our paper is the following:

\begin{theorem} \label{The bound}
Assume that conditions $({\bf H_{\eta}})$, $({\bf H1})$-$({\bf H3})$ hold, and the coefficients $\sigma$, $b$ are $\mathcal{C}^{\infty}$ functions with bounded derivatives of all orders.
Then for all $(t,x)\in(0, T]\times \RR^d$, the law of the random vector $u(t,x)$ has a smooth density, denoted by $p_{t,x}(y)$, which satisfies that, for all $y\in\RR^m$,
\begin{equation}
C_1\Phi(t)^{-\frac{m}{2}}\exp\left(-\frac{|y|^2}{C_2\Phi(t)}\right)\leq p_{t,x}(y)\leq C_3\Phi(t)^{-\frac{m}{2}}\exp\left(-\frac{(|y|-C_4 T)^2}{C_5\Phi(t)}\right), \label{lower and upper}
\end{equation}
where $\Phi(t)=\int^{t}_{0}\int_{\RR^d}|\cF\Gamma(s)(\xi)|^2\mu(d\xi)ds$, $C_i, i=1,\ldots, 5$ are positive constants that only depend on $T, \sigma$ and $b$.
\end{theorem}

The proof of this theorem will be finished by the following two subsections. One studies the lower bound and the other studies the upper bound.


\subsection{\bf The lower bound}

The concept of uniformly elliptic random vector was used to obtain the lower bound for the density of a random vector (see \cite{K} or \cite{NuaEQ}). If the solution $u(t,x)$ of Eq.\eref{Eq} is a uniformly elliptic $m$-dimensional random vector, the low bound of the density of $u(t,x)$ will be got by \cite[Theorem 2.3]{NuaEQ}. Now, we recall the definition of uniformly elliptic random vector (see \cite[Definition 2.2]{NuaEQ}).

\begin{definition}\label{ue}
Let $F$ be a non-degenerate $m$-dimensional $\cF_t$-measurable random vector. $F$ is called uniformly elliptic if there exists an $\varepsilon>0$ such that for any partition $\pi_N=\{0= t_0< t_1< \cdots< t_N=t\}$ whose norm $\|\pi_N\|:=\max\{|t_{i+1}-t_{i}|; i=0, \ldots, N-1\}$ is smaller than $\varepsilon>0$ and $\|\pi_N\|\to0$ as $N\to\infty$, there exists a sequence of smooth random vectors
$(F_n)_{n=0,\ldots,N}$ such that $F_N=F$, $\cF_{t_n}$-measurable $F_n$ belongs to $(\overline{\mathbb{D}}^{\infty}_{t_{n-1}, t_n})^m$ and $F_n$ can be written in the following form:
\begin{equation}
F_n=F_{n-1}+I_n(h)+G_n,~~~~n=1,\ldots, N,
\end{equation}
where the random vectors $I_n(h)$ and $G_n$ satisfy the following conditions:

(A1) $G_n$ is an $\cF_{t_n}$-measurable and belongs to $(\overline{\mathbb{D}}^{\infty}_{t_{n-1}, t_n})^m$, and there exists an element $g\in \mathcal{H}_T$ with $\|g(s)\|_{\mathcal{H}}>0$ (a.s. $s$) such that,
for all $k\in \NN$ and $p\geq 1$,
\begin{eqnarray}
\|G_n\|^{t_{n-1}, t_n}_{k, p}\leq C \Delta_{n-1}(g)^{1/2+\gamma}\quad a.s.,\label{G}
\end{eqnarray}
for some $\gamma>0$, where
\[
0<\Delta_{n-1}(g):=\int^{t_n}_{t_{n-1}}\|g(s)\|^2_{\mathcal{H}}ds< \infty,\quad n=1,\ldots, N.
\]

(A2) Random vector $I_n(h)$ with the component:
\[
I^{i}_n(h)=\int^{t_n}_{t_{n-1}}\int_{\RR^d}h_{i}(s, y) W(ds, dy),\quad i=1,\ldots, m,
\]
where $h_i$ is a smooth $\cF_{t_{n-1}}$-predictable $\mathcal{H}^q_{t_{n-1}, t_n}$-valued process.
For $k\in\NN$, $p\geq 1$ and $i\le m$, there exists a constant $C>0$ such that
\[
\|F^{i}_n\|_{k,p}+\sup_{\omega\in\Omega}\|h_{i}\|_{t_{n-1}, t_n}(\omega)\leq C.
\]

(A3) Let $A=(a_{i,j})$ denote the $m\times m$ matrix defined by
\[
a_{i,j}=\Delta_{n-1}(g)^{-1}\int^{t_n}_{t_{n-1}}\langle h_i(s), h_j(s)\rangle_{\mathcal{H}^q} ds.
\]
There exist strictly positive constants $C_1$ and $C_2$ such that, for all $\xi\in\RR^m$,
\[
C_1|\xi|^2\leq \xi^{T}A\xi\leq C_2|\xi|^2.
\]

(A4) There is a constant $C$ such that, for $p> 1$ and $\rho\in (0,1]$,
\[
\EE_{t_{n-1}} \left[\text{det}(M^{t_{n-1}, t_n}_{I_n(h)+\rho G_n})^{-p}\right]\leq C\Delta_{n-1}(g)^{-mp}\quad a.s..
\]
\end{definition}

The following theorem shows the lower bound in Theorem \ref{The bound}.

\begin{theorem}
Under the assumptions in Theorem \ref{The bound}. Then the solution $u(t,x)$ of Eq.\eref{Eq} is an $m$-dimensional uniformly elliptic random vector. And the density $p_{t,x}(y)$ of $u(t,x)$ satisfies:
\begin{equation}
p_{t,x}(y)\geq C_1\Phi(t)^{-\frac{m}{2}}\exp\left(-\frac{|y|^2}{C_2\Phi(t)}\right),  \quad \forall y\in\RR^m\label{lower},
\end{equation}
where $\Phi(t)=\int^{t}_{0}\int_{\RR^d}|\cF\Gamma(s)(\xi)|^2\mu(d\xi)ds$.
\end{theorem}

Proof. Refer to \cite[Theorem 2.3]{NuaEQ}, it suffices to check that $u(t,x)$ is a $m$-dimensional uniformly elliptic random vector with $g(\cdot):=\Gamma(t-\cdot)$ in (A1).

We consider a partition $0=t_0< t_1\cdots < t_N=t$ with $\sup_{1\le i\le N}(t_i-t_{i-1})\to 0$ as $N\to\infty$, and define, for $ i=1,\ldots, m$,
\begin{eqnarray*}
F^{i}_n&=&\int^{t_n}_{0}\int_{\RR^d}\sum^q_{j=1}\Gamma(t-s, x-y)\sigma_{ij}(u(s,y))W^{j}(ds,dy)\\
&&+\int^{t_n}_{0}\int_{\RR^d}b_i(u(s, y))\Gamma(t-s, x-y)dyds.
\end{eqnarray*}
It is obvious that $F_n:=(F^1_n,\ldots, F^m_n)$ is $\cF_{t_n}$-measurable ($n=0,\ldots, N$), $F_0=0$ and $F_N=u(t,x)$. Moreover, $F_n\in (\mathbb{D}^{\infty})^{m}$ and, for all $k\in \NN$ and $p>1$, the norm $\|F^i_n\|_{k, p}$ can be uniformly bounded with respect to $(t, x)\in (0, T]\times \RR^d$.

Let $g(s)=\Gamma(t-s)$. \eref{2.2} and \eref{ie.4.1} imply
\begin{eqnarray}
0<\Delta_{n-1}(g):=\int^{t_{n}}_{t_{n-1}}\|g(s)\|^2_{\mathcal{H}}ds=\int^{t_{n}}_{t_{n-1}}\int_{\RR^d}|\cF\Gamma(t-s)(\xi)|^2 \mu(d\xi)ds<\infty.\label{g}
\end{eqnarray}

Next, we intend to decompose $F_n$ in the form (3.6). For any $i=1,\ldots, m$, we have
\begin{eqnarray*}
F^{i}_n-F^{i}_{n-1}&=&\int^{t_n}_{t_{n-1}}\int_{\RR^d}\sum^q_{j=1}\Gamma(t-s, x-y)\sigma_{ij}(u(s,y))W^j(ds,dy)\\
&&+\int^{t_n}_{t_{n-1}}\int_{\RR^d}b_i(u(s, y))\Gamma(t-s,x-y)dyds\\
&=&\int^{t_n}_{t_{n-1}}\int_{\RR^d}\sum^q_{j=1}\Gamma(t-s, x-y)\sigma_{ij}(u_{n-1}(s,y))W^j(ds,dy)\\
&&+\int^{t_n}_{t_{n-1}}\int_{\RR^d}b_i(u(s, y))\Gamma(t-s, x-y)ds\\
&&+\int^{t_n}_{t_{n-1}}\int_{\RR^d}\sum^q_{j=1}\Gamma(t-s, x-y)\left[\sigma_{ij}(u(s,y))-\sigma_{ij}(u_{n-1}(s,y))\right]W^j(ds,dy),
\end{eqnarray*}
where $u_{n-1}(s,y)=(u^1_{n-1}(s,y),\ldots,u^{m}_{n-1}(s,y))$ is defined by
\begin{eqnarray*}
u^{i}_{n-1}(s,y)&=&\int^{t_{n-1}}_{0}\int_{\RR^d}\sum^{q}_{j=1}\Gamma(s-r, y-z)\sigma_{ij}(u(r,z))W^j(dr,dz)\\
&&+\int^{t_{n-1}}_{0}\int_{\RR^d}b_i(u(r,z))\Gamma(s-r,y-z)dz dr,\ \ i=1,\ldots, m,
\end{eqnarray*}
for $(s, y)\in [t_{n-1}, t_n]\times \RR^d$. Clearly, $u^i_{n-1}(s, y)$ is $\cF_{t_{n-1}}$-measurable and belongs to $\mathbb{D}^{\infty}$.

Hence, we can obtain a decomposition of $F_n$:
\[
F_n=F_{n-1}+I_n(h)+G_n,
\]
where $I_n(h):=(I^{1}_n(h),\ldots, I^{m}_n(h))$ with
\[
I^{i}_n(h):=\int^{t_n}_{t_{n-1}}\int_{\RR^d}\sum^q_{j=1}h_{ij}(s,y) W^j(ds, dy)
\]
and
$$h_{ij}(s,y):=\Gamma(t-s, x-y)\sigma_{ij}(u_{n-1}(s,y)),$$
and $G_n:=(G^{1}_n,\ldots, G^{m}_n)$ with
\begin{eqnarray*}
G^{i}_n:=\!\!\!\!\!\!\!\!&&\int^{t_n}_{t_{n-1}}\int_{\RR^d}b_i(u(s, y))\Gamma(t-s, x-y)ds\\
&&+\int^{t_n}_{t_{n-1}}\int_{\RR^d}\sum^q_{j=1}\Gamma(t-s, x-y)[\sigma_{ij}(u(s,y))-\sigma_{ij}(u_{n-1}(s,y))]W^j(ds,dy).
\end{eqnarray*}

Firstly, \eref{G} is satisfied by \cite[Lemma 4.1]{NuaEQ}. This and \eref{g} yield that (A1) holds.
Secondly, the boundedness of $\sup_{\omega\in \Omega}\|h_i\|_{t_{n-1},t_n}$ is a consequence of the condition (A3) (see step 1 below).
All conditions in (A2) are fulfilled by the boundedness of $\sup_{\omega\in \Omega}\|h_i\|_{t_{n-1},t_n}$ and $F^i_n\in \mathbb{D}^{\infty}$.
The remaining is to check the conditions in (A3) and (A4), which will be done by the following two steps.

{\it Step 1}. We will prove that there exist two positive constants $C_1$ and $C_2$ such that, for all $\xi\in\RR^m$,
\[
C_1\xi^{T}\xi\leq\xi^{T}A\xi\leq C_2\xi^{T}\xi,
\]
where $A:=(a_{i,j})$ is the $m\times m$ matrix defined by
\[
a_{i,j}:=\Delta_{n-1}(g)^{-1}\int^{t_n}_{t_{n-1}}\langle h_i(s), h_j(s)\rangle_{\mathcal{H}^q} ds.
\]
Without loss of generality, we assume $|\xi|=1$. By \eref{lower of sigma}, we have
\begin{eqnarray*}
\xi^{T}A\xi&=&\Delta_{n-1}(g)^{-1}\int^{t_n}_{t_{n-1}}\left\|\sum^m_{i=1}h_i(s)\xi_i\right\|^2_{\mathcal{H}^q}ds\\
&\geq&\Delta_{n-1}(g)^{-1}\int^{t_n}_{t_{n-1}}\|\Gamma(t-s, x-\ast)\|^2_{\mathcal{H}}\inf_{x,y\in\RR^d}\sum^m_{i,j=1}\sum^q_{k=1}\sigma_{ik}(x)\sigma_{jk}(y)\xi_i\xi_jds\\
&\geq&C_1\Delta_{n-1}(g)^{-1}\int^{t_n}_{t_{n-1}}\|\Gamma(t-s, x-\ast)\|^2_{\mathcal{H}}ds=C_1.
\end{eqnarray*}
Meanwhile, \eref{upper of sigma} yields
\begin{eqnarray*}
\xi^{T}A\xi&=&\Delta_{n-1}(g)^{-1}\int^{t_n}_{t_{n-1}}\left\|\sum^m_{i=1}h_i(s)\xi_i\right\|^2_{\mathcal{H}^q}ds\\
&\leq&\Delta_{n-1}(g)^{-1}\int^{t_n}_{t_{n-1}}\|\Gamma(t-s, x-\ast)\|^2_{\mathcal{H}}\sup_{x, y\in\RR^d}\sum^m_{i,j=1}\sum^q_{k=1}\sigma_{ik}(x)\sigma_{jk}(y)\xi_i\xi_jds\\
&\leq&C_2\Delta_{n-1}(g)^{-1}\int^{t_n}_{t_{n-1}}\|\Gamma(t-s, x-\ast)\|^2_{\mathcal{H}}ds=C_2.
\end{eqnarray*}
The condition (A3) is satisfied.

{\it Step 2}.  We check the condition (A4), i.e., for any $p>0$, there exists a constant $C>0$ such that
\begin{equation}
\EE_{t_{n-1}} \left[\text{det}(M^{t_{n-1}, t_n}_{I_n(h)+\rho G_n})^{-p}\right]\leq C\Delta_{n-1}(g)^{-mp}\quad a.s..\label{M_t}
\end{equation}
In fact, by \cite[Lemma 2.3.1]{NuaD}, it is sufficient to prove that for any $q\geq2$, there exists $\varepsilon_0=\varepsilon_0(q)>0$ such that, for all $\varepsilon\leq \varepsilon_0$,
\[
\sup_{|\xi|=1}\mathbb{P}_{t_{n-1}}\left\{\xi^{T}(\Delta_{n-1}(g)^{-1}M^{t_{n-1}, t_n}_{I_n(h)+\rho G_n})\xi\leq \varepsilon\right\}\leq \varepsilon^q, \quad a.s..
\]
The term $I_n(h)+\rho G_n$ can be split as follows:
\[
I_n(h)+\rho G_n=\rho(I_n(h)+G_n)+(1-\rho)I_n(h)=\rho(F_n-F_{n-1})+(1-\rho)I_n(h).
\]
Thus, for any $r\in[t_{n-1}, t_n]$, we have
\begin{eqnarray*}
D_{r,\ast}(I^i_n(h)+\rho G^i_n)
=\!\!\!\!\!\!\!\!&& \rho\bigg\{\Gamma(t-r, x-\ast)\sigma_{i}(u(r, \ast))\\
&&+\int^{t_n}_r\!\!\!\!\int_{\RR^d}\sum^{q}_{l=1}\Gamma(t-s, x-y)D_{r, \ast}(\sigma_{il}(u(s,y)))W^l(ds,dy)\\
&&+\int^{t_n}_r\!\!\!\!\int_{\RR^d}\Gamma(t-s, x-y)D_{r, \ast}(b_i(u(s,y)))dyds\bigg\}\\
&&+(1-\rho)\bigg\{\Gamma(t-r, x-\ast)\sigma_{i}(u_{n-1}(r, \ast))\\
&&~~~~+\int^{t_n}_r\!\!\!\!\int_{\RR^d}\sum^q_{l=1}\Gamma(t-s, x-y)D_{r, \ast}(\sigma_{il}(u_{n-1}(s,y)))W^l(ds,dy)\bigg\}\\
=\!\!\!\!\!\!\!\!&& \rho\Gamma(t-r, x-\ast)\sigma_{i}(u(r, \ast))+(1-\rho)\Gamma(t-r, x-\ast)\sigma_{i}(u_{n-1}(r, \ast))\\
&&+\rho\bigg\{\int^{t_n}_r\!\!\!\!\int_{\RR^d}\sum^{q}_{l=1}\Gamma(t-s, x-y)D_{r, \ast}(\sigma_{il}(u(s,y)))W^l(ds,dy)\\
&&~~~~+\int^{t_n}_r\!\!\!\!\int_{\RR^d}\Gamma(t-s, x-y)D_{r, \ast}(b_i(u(s,y)))dyds\bigg\},
\end{eqnarray*}
the last equality comes from $D_{r}(u_{n-1}(s,y)))=0$ for $r\in(t_{n-1}, t_n]$  when $u_{n-1}(s,y)$ is $\mathcal{F}_{t_{n-1}}$-measurable.
Therefore, for $\delta\in (0, t_n-t_{n-1}]$, we have
\begin{eqnarray*}
\xi^{T}(\Delta_{n-1}(g)^{-1}M^{t_{n-1}, t_n}_{I_n(h)+\rho G_n})\xi
=\!\!\!\!\!\!\!\!&&\int^{t_n}_{t_{n-1}}\Delta_{n-1}(g)^{-1}\left\|\sum^m_{i=1}D_{r,\ast}(I^i_n(h)+\rho G^i_n)\xi_i\right\|^2_{\mathcal{H}^q}dr\\
\geq\!\!\!\!\!\!\!\!&&\int^{t_n}_{t_{n}-\delta}\Delta_{n-1}(g)^{-1}\left\|\sum^m_{i=1}D_{r,\ast}(I^i_n(h)+\rho G^i_n)\xi_i\right\|^2_{\mathcal{H}^q}dr\\
\geq\!\!\!\!\!\!\!\!&& \frac{1}{2}\Delta_{n-1}(g)^{-1}\mathcal{B}_1- \Delta_{n-1}(g)^{-1}\mathcal{B}_2,
\end{eqnarray*}
where
\begin{eqnarray*}
\mathcal{B}_1:=\!\!\!\!\!\!\!\!&&\int^{t_n}_{t_{n}-\delta}\left\|\sum^m_{i=1}\left[\rho\Gamma(t-r, x-\ast)\sigma_i(u(r,\ast))+(1-\rho)\Gamma(t-r, x-\ast)\sigma_i(u_{n-1}(r,\ast))\right]
\xi_i\right\|^2_{\mathcal{H}^q}dr\,,
\end{eqnarray*}
\begin{eqnarray*}
\mathcal{B}_2:=\!\!\!\!\!\!\!\!&&\int^{t_n}_{t_{n}-\delta}\|a(r, t, x, \ast)\|^2_{\mathcal{H}^q}dr
\end{eqnarray*}
and
\begin{eqnarray*}
a(r, t, x, \ast)=\!\!\!\!\!\!\!\!&&\sum^{m}_{i=1}\int^{t_n}_r\int_{\RR^d}\rho\sum^{q}_{l=1}\Gamma(t-s, x-y)D_{r, \ast}(\sigma_{il}(u(s,y)))W^l(ds,dy)\xi_i\\
&&+\sum^{m}_{i=1}\int^{t_n}_r\int_{\RR^d}\rho\Gamma(t-s, x-y)D_{r, \ast}(b_i(u(s,y)))dyds\xi_i\\
:=\!\!\!\!\!\!\!\!&& I_1+I_2.
\end{eqnarray*}

Put $I_0(\delta):=\int^{t_n}_{t_n-\delta}\int_{\RR^d}|\Gamma(t-r)(\xi)|^2\mu(d\xi) dr$.
By (\ref{lower of sigma}), we have
\begin{eqnarray*}
&&\bigg\|\sum^m_{i=1}\big[\rho\Gamma(t-r, x-\ast)\sigma_i(u(r,\ast))+(1-\rho)\Gamma(t-r, x-\ast)\sigma_i(u_{n-1}(r,\ast))\big]\xi_i\bigg\|^2_{\mathcal{H}^q}\,, \\
=\!\!\!\!\!\!\!\!&&\!\!\int_{\RR^d}\!\!\int_{\RR^d}\!\sum^m_{i,j=1}\sum^q_{k=1}\big[\rho\Gamma(t\!-\!r, x\!-\!y)\sigma_{ik}(u(r,y))\!+\!(1\!-\!\rho)\Gamma(t\!-\!r, x\!-\!y)\sigma_{ik}(u_{n-1}(r,y))\big]f(y\!-\!z)\\
&&\times\big[\rho\Gamma(t-r, x-z)\sigma_{jk}(u(r,z))+(1-\rho)\Gamma(t-r, x-z)\sigma_{jk}(u_{n-1}(r,z))\big]\xi_i\xi_j dy dz\\
\geq\!\!\!\!\!\!\!\!&&[\rho^2+2\rho(1-\rho)+(1-\rho)^2]\|\Gamma(t-r, x-\ast)\|^2_{\mathcal{H}}\inf_{x,y\in\RR^d}\sum^m_{i,j=1}\sum^q_{k=1}\sigma_{ik}(x)\sigma_{jk}(y)\xi_i\xi_j\\
\geq\!\!\!\!\!\!\!\!&& C_1\|\Gamma(t-r, x-\ast)\|^2_{\mathcal{H}}.
\end{eqnarray*}
This implies that $\mathcal{B}_{1}\geq C_1 I_0(\delta)$.

Next, we are going to estimate the terms $\EE_{t_{n-1}}\|I_i\|^{2p}_{t_n-\delta, t_n}$ ($i=1,2$) for $p>1$.

By the boundedness of the partial derivative of $\sigma_{ij}$, BDG inequality and H\"{o}lder's inequality, we have
\begin{eqnarray}
\EE_{t_{n-1}}\|I_1\|^{2p}_{t_n-\delta, t_n}
\leq\!\!\!\!\!\!\!\!&& C \sum^{m}_{i=1}\EE_{t_{n-1}}\!\left\|\int^{t_n}_{t_n-\delta}\!\int_{\RR^d}\sum^{q}_{l=1}\Gamma(t-s, x-y)D(\sigma_{il}(u(s,y)))W^l(ds,dy)\right\|^{2p}_{t_n-\delta, t_n}\nonumber\\
\leq\!\!\!\!\!\!\!\!&& C \sum^{m}_{i=1}I^{p-1}_0(\delta)\int^{t_n}_{t_n-\delta}\sup_{z\in\RR^d}\EE_{t_{n-1}}(\|Du_i(s,z)\|^{2p}_{t_n-\delta, t_n})\|\Gamma(t-s, \ast)\|^2_{\mathcal{H}}ds \nonumber\\
\leq\!\!\!\!\!\!\!\!&& C \sum^{m}_{i=1}I^{p}_0(\delta)\sup_{s\in[t_n-\delta, t_n],z\in\RR^d}\EE_{t_{n-1}}(\|Du_i(s,z)\|^{2p}_{t_n-\delta, t_n})\nonumber\\
\leq\!\!\!\!\!\!\!\!&& C I^{p}_0(\delta)\Phi^p(\delta),\quad a.s.,\label{I_1}
\end{eqnarray}
where the last inequality comes from  \eref{2.7}.

Similarly, the boundedness of the partial derivative of $b_{i}$ implies
\begin{eqnarray}
\EE_{t_{n-1}}\|I_2\|^{2p}_{t_n-\delta, t_n}
\leq\!\!\!\!\!\!\!\!&& C \sum^{m}_{i=1}\EE_{t_{n-1}}\left\|\int^{t_n}_{t_n-\delta}\int_{\RR^d}\Gamma(t-s, x-y)D(b_i(u(s,y)))dyds\right\|^{2p}_{t_n-\delta, t_n}\nonumber\\
\leq\!\!\!\!\!\!\!\!&& C \sum^{m}_{i=1}\overline{I}^{p-1}_0(\delta)\int^{t_n}_{t_n-\delta}\sup_{z\in\RR^d}\EE_{t_{n-1}}(\|Du_i(s,z)\|^{2p}_{t_n-\delta, t_n})\int_{\RR^d}\Gamma(t-s, z)dzds \nonumber\\
\leq\!\!\!\!\!\!\!\!&& C \sum^{m}_{i=1}\overline{I}^{p}_0(\delta)\sup_{s\in[t_n-\delta, t_n],z\in\RR^d}\EE_{t_{n-1}}(\|Du_i(s,z)\|^{2p}_{t_n-\delta, t_n})\nonumber\\
\leq\!\!\!\!\!\!\!\!&& C \overline{I}^{p}_0(\delta)\Phi^p(\delta), \quad a.s.,\label{I_2}
\end{eqnarray}
where
\begin{eqnarray}
\overline{I}_0(\delta):=\int^{t_n}_{t_n-\delta}\int_{\RR^d}\Gamma(t-s, z)dzds\leq C\delta. \label{3.12}
\end{eqnarray}
\eref{I_1} and \eref{I_2} show
\begin{eqnarray}
\EE_{t_{n-1}}|\mathcal{B}_2|^p\leq C \Phi(\delta)^p(I_0(\delta)^p+\overline{I}_0(\delta)^p),\quad a.s..\label{ie.4.5}
\end{eqnarray}
Hence, by \eref{ie.4.5} and the (conditional) Chebyshev's inequality, we obtain
\begin{eqnarray*}
&&\sup_{|\xi|=1}\mathbb{P}_{t_{n-1}}\left\{\xi^{T}(\Delta_{n-1}(g)^{-1}M^{t_{n-1}, t_n}_{I_n(h)+\rho G_n})\xi\leq \varepsilon\right\}\\
\leq\!\!\!\!\!\!\!\!&&\mathbb{P}_{t_{n-1}}\left\{\mathcal{B}_2\geq \frac{C_1}{2}I_0(\delta)-\Delta_{n-1}(g)\varepsilon\right\}\\
\leq\!\!\!\!\!\!\!\!&& C\left(\frac{C_1 I_0(\delta)}{2\Delta_{n-1}(g)}-\varepsilon\right)^{-p}(\Delta_{n-1}(g))^{-p}\left[\Phi(\delta)^p(I_0(\delta)^p+\overline{I}_0(\delta)^p)\right].
\end{eqnarray*}

Now, taking a small enough $\varepsilon_0$ if necessary, we choose $\delta=\delta(\varepsilon)$ such that $\frac{C_1 I_0(\delta)}{4\Delta_{n-1}(g)}=\varepsilon$.
By \eref{ie.4.1}, \eref{ie.4.2} and \eref{3.12}, we have
\begin{eqnarray*}
\sup_{|\xi|=1}\mathbb{P}_{t_{n-1}}\left\{\xi^{T}(\Delta_{n-1}(g)^{-1}M^{t_{n-1}, t_n}_{I_n(h)+\rho G_n})\xi\leq \varepsilon\right\}
\leq\!\!\!\!\!\!\!\!&& CI_0(\delta)^{-p}\left[\Phi(\delta)^p(I_0(\delta)^p+\overline{I}_0(\delta)^p)\right]\\
\leq\!\!\!\!\!\!\!\!&& C\delta^{(1-\eta)p}\leq C\varepsilon^{(1-\eta)p},
\end{eqnarray*}
We use the fact:
\[
\delta\leq C I_0(\delta)\leq C\Delta_{n-1}(g)\varepsilon\leq C\varepsilon
\]
in the last inequality.
\rule{0.5em}{0.5em}


\subsection{\bf The upper bound}

This subsection is devoted to prove the upper bound of the joint density. To do this, we will use a classical method based on the density formula provided by the integration by parts formula of the Malliavin calculus (see \cite[Corollary 3.2]{GMEN} or \cite[Proposition 2.1.5]{NuaD}).

We first consider the continuous $\RR^m$-valued martingale $\{Z_a, \cF_a, 0\leq a\leq t\}$ defined by
\[
Z^i_a:=\int^{a}_{0}\int_{\RR^d}\sum^{q}_{j=1}\Gamma(t-s, x-y)\sigma_{ij}(u(s,y))W^j(ds, dy),\quad i=1,\ldots, m.
\]
Notice that
\[
\langle Z\rangle_t=\sum^{m}_{i=1}\|\Gamma(t-\cdot, x-\ast)\sigma_{i}(u(\cdot,\ast))\|^2_{\mathcal{H}^q_t}.
\]
By \eref{upper of sigma}, there exists some positive constant $C_1$ depending on $\sigma$ and $m$ such that $\langle Z\rangle_t\leq C_1\Phi(t)$.

Since $\Gamma(t,x)$ is a Gaussian density and $b_i$ ($i=1,\ldots, m$) are bounded, we have that, for $t\in(0, T]$ and $x\in\RR^d$,
\begin{equation}
\left|\int^t_0\int_{\RR^d}\Gamma(t-s, x-y)b_i(u(s,y))dy ds\right|\leq C T,\label{4.9}
\end{equation}
where $C$ is a constant depending on $b$.

Next, we consider the expression of the joint density of a non-degenerate random vector.
Using the integration by part formula of the Malliavin calculus (see \cite[Corollary 3.2]{GMEN}), we have the following expression of the joint density $p_{t,x}(\cdot)$ of $u(t,x)$,
\begin{equation}
p_{t, x}(y)=(-1)^{m-\text{card}(\SS)}\EE\left[{\bf 1}_{\{u_i(t,x)>y_i,~i\in\SS;~u_i(t,x)<y_i,~i\not\in\SS;~i=1,\ldots,m\}}H_{(1,2,\ldots, m)}(u(t,x),1)\right], \quad y\in\RR^d, \nonumber
\end{equation}
where $\SS$ be a subset of $\{1,\ldots, m\}$, $\text{card}(\SS)$ denotes the cardinality of $\SS$, the random variables $H_{\alpha}(F,G)$ are recursively given by
\[
H_{(i)}(F, G):=\sum^{m}_{j=1}\delta(G(M^{-1}_F)_{ij}DF^j),
\]
\[
H_{\alpha}(F, G):= H_{(\alpha_m)}(F, H_{(\alpha_1,\alpha_2,\ldots, \alpha_{m-1})}(F,G))
\]
for any $F\in (\mathbb{D}^{\infty})^m$, $G\in\mathbb{D}^{\infty}$ and $\alpha=(\alpha_1, \ldots, \alpha_m)\in\{1,\ldots, m\}^m$.
Then, by H\"{o}lder's inequality and \eref{4.9}, we get
\begin{eqnarray}
p_{t, x}(y)
\leq\!\!\!\!\!\!\!\!&&\PP\left\{|u(t, x)|>|y|\right\}^{1/2}\left\{\EE[H_{(1,2,\ldots, m)}(u(t,x),1)]^2\right\}^{1/2}\nonumber\label{iq.4.11}\\
\leq\!\!\!\!\!\!\!\!&&\PP\{|Z_t|>|y|-CT\}^{1/2}\left\{\EE[H_{(1,2,\ldots, m)}(u(t,x),1)]^2\right\}^{1/2}.
\end{eqnarray}
$\langle Z\rangle_t\leq C_1\Phi(t)$ and the exponential martingale inequality (see for instance \cite[Section A2]{NuaD}) imply
\begin{eqnarray}
\PP\{|Z_t|>|y|-CT\}\leq 2\exp\left\{-\frac{(|y|-CT)^2}{2C_1\Phi(t)}\right\}.\label{iq.4.12}
\end{eqnarray}
Meanwhile, by \cite[Lemma 3.4]{NuaEQ} and \eref{M_t}, we have the following two estimates:

(i) $\left\|\mathbb{D}^k(u_i(t,x))\right\|_{L^p(\Omega, \mathcal{H}^k_t)}\leq C\Phi(t)^{1/2}$,\quad i=1,\ldots, m,

(ii)$\left\|\text{det}(M_{u(t,x)})^{-1}\right\|_{L^p(\Omega)}\leq C\Phi(t)^{-m}$.\\
Then \cite[Proposition 3.3]{GMEN} implies that there exists a constant $C>0$ such that
\begin{equation}
\|H_{(1,2,\ldots, m)}(u(t,x),1)\|_{L^p(\omega)}\leq C\Phi(t)^{-\frac{m}{2}}.\label{iq.4.13}
\end{equation}
Hence, \eref{iq.4.11}-\eref{iq.4.13} yield the desired upper bound in \eref{lower and upper}.
\rule{0.5em}{0.5em}

\section{Examples}

Let $\Gamma(r,x)=(2\pi r)^{-d/2}e^{-\frac{|x|^2}{2r}}$ be the fundamental solution for the heat equation on $\RR^d$.
We will give some examples of  covariance functions $f$ satisfying hypotheses ${(\bf H_{\eta})}$, ${(\bf H1)}$ and ${(\bf H2)}$.

\vspace{0.2cm}
\medskip
\noindent
\textit {\textbf{Riesz kernel}}. Let  $f(x)=|x|^{-\gamma}$ with $0 < \gamma < 2 \wedge d$ and $\mu(d\xi)=C_{d,\gamma}|\xi|^{\gamma-d}d\xi$.
Then ${(\bf H_{\eta})}$ holds for any $\eta>\frac{\gamma}{2}$. \eref{gamma_1} and \eref{gamma_2} are satisfied with $0< \gamma_1<\frac{2-\gamma}{4}$
and $0< \gamma_2 < \frac{2-\gamma}{2}$ respectively.  According to \cite{NuaE}, $({\bf H1})$  holds with $\beta = \frac{2-\gamma}{2}$. (${\bf H 2}$)
holds with  $\beta_1= \frac {2-\gamma} 2 +\frac{\gamma_2}{2}$, $\beta_2=\frac {2-\gamma} 2 +\gamma_1$.

\vspace{0.2cm}
\medskip
\noindent
\textit {\textbf{Bessel kernel}}.
Let $f(x)=\int^{\infty}_{0}u^{\frac{\alpha-d-2}{2}}e^{-u}e^{-\frac{|x|^2}{4u}}du$ for $d-2<\alpha<d$ and $\mu(d\xi)= c_{\alpha,d} (1+ |\xi|^2) ^{-\frac {\alpha}2 } d\xi$.
Then ${(\bf H_{\eta})}$ holds for $\eta>\frac{d-\alpha}{2}$. According to \cite{NuaE}, \eref{gamma_1} and \eref{gamma_2} are satisfied with $0< \gamma_1<\frac{2-d+\alpha}{4}$
and $0< \gamma_2 < \frac{2-\gamma+\alpha}{2}$ respectively.  For $\varepsilon <1$, we have
\begin{eqnarray*}
\int_0^{\varepsilon} \int_{\RR^d} |\mathcal{F}\Gamma(r)(\xi)|^2\mu(d\xi)dr
=\!\!\!\!\!\!\!\!&&C \int_0^{\varepsilon} \int_{\RR^d}e^{-r|\xi|^2}  (1+ |\xi|^2)  ^{-\frac   {\alpha}{2}} d\xi dr \\
=\!\!\!\!\!\!\!\!&&C \int_0^{\varepsilon}\int_{\RR^d}e^{-|\theta|^2}\frac{r^{\frac{\alpha-d}{2}}}{(|\theta|^2+ r)^{ \frac {\alpha}2 }}d\theta dr\\
\geq\!\!\!\!\!\!\!\!&& C \int_0^{\varepsilon}r^{\frac{\alpha-d}{2}}dr\int_{\RR^d}e^{-|\theta|^2}\frac{1}{(|\theta|^2+1)^{ \frac {\alpha}2}}d\theta\\
=\!\!\!\!\!\!\!\!&& C \varepsilon^{\frac{\alpha-d}{2}+1}\,.
\end{eqnarray*}
This yields that $({\bf H1})$ is satisfied with $\beta = \frac{\alpha-d}{2}+1$. To show (${\bf H2}$), Using $f(x)\leq C |x|^{-d+\alpha}$
for $x \in \RR^d$ (see \cite[Proposition 6.1.5]{Gra}) and proceeding as in the case of the Riesz kernel with $\beta= d-\alpha$,
we obtain that (${\bf H2}$) holds with $\beta_1 = \frac {2+\alpha -d} 2 +\frac{\gamma_2}{2}$ and $\beta_2=\frac {2+\alpha -d} 2 +\gamma_1$.

\vspace{0.2cm}
\medskip
\noindent
{\textit{\textbf{Fractional kernel}}.
Let $f(x)=\prod_{j=1}^d |x_j|^{2H_j-2}$ with $\sum_{j=1}^d H_j > d-1$ for $ \frac{1}{2} < H_j < 1$ ($1\leq j \leq d$). It is clear that all of our theory still works for this case although
$f(x)$ is continuous on $\RR^d \setminus \{0\}$ only. Then we have $\mu(d\xi)=C_H\prod_{j=1}^d |\xi_j|^{1-2H_j}d\xi$, where $C_H$ only depends on $H:=(H_1,H_2,\dots,H_d)$.

${(\bf H_{\eta})}$ holds for $\eta>d-\sum_{j=1}^d H_j$. According to \cite{NuaE}, \eref{gamma_1} and \eref{gamma_2} are satisfied with
$0< \gamma_1<(1/2)(\sum_{j=1}^d H_j -d+1)$ and $0< \gamma_2 < \sum_{j=1}^d H_j -d+1$ respectively.  Using the change of variable $\sqrt{t}\xi\to\xi$, we obtain
 \begin{eqnarray*}
 \int_{0}^\varepsilon \int_{\RR^d} |\mathcal{F}\Gamma(t)(\xi)|^2\mu(d\xi)dt
 =\int_0^\varepsilon \int_{\RR^d} e^{-t|\xi|^2}\prod_{j=1}^d |\xi_j|^{1-2H_j}d\xi dt =C \varepsilon^{\sum_{j=1}^d H_j -d+1}.
  \end{eqnarray*}
Thus,  $({\bf H1})$ is satisfied with $\beta=\sum_{j=1}^d H_j -d+1$. From $|x|^{\alpha}\Gamma(r,x)\leq C r^{\frac{\alpha}{2}} \Gamma(2r,x)$,   we have that, for all $x\in \RR^d$,
\begin{eqnarray*}
  \int_0^{\varepsilon} \langle |*|^{\gamma_2} \Gamma(r,*), \Gamma(r,*) \rangle_{\mathcal{H}}dr
=\!\!\!\!\!\!\!\!&& \int_0^{\varepsilon} \int_{\RR^d}\int_{\RR^d} |x|^{\gamma_2} \Gamma(r,x) \Gamma(r,y)  \prod_{j=1}^d |x_j-y_j|^{2H_j-2}dx dy dr\\
\leq\!\!\!\!\!\!\!\!&&C \int_0^{\varepsilon} \int_{\RR^d}\int_{\RR^d} r^{\frac{\gamma_2}{2}}  \Gamma(2r,x) \Gamma(r,y) \prod_{j=1}^d |x_j-y_j|^{2H_j-2}dx dy dr\\
=\!\!\!\!\!\!\!\!&&C \int_0^{\varepsilon} \int_{\RR^d} r^{\frac{\gamma_2}{2}} e^{-\frac{3r}{2} |\xi|^2} \prod_{j=1}^d |\xi_j|^{1-2H_j} d\xi dr\\
=\!\!\!\!\!\!\!\!&&C \int_0^{\varepsilon} r^{\frac{\kappa_2}{2}+\sum_{j=1}^d H_j -d}dr= C \varepsilon^{\frac{\gamma_2}{2}+\sum_{j=1}^d H_j -d+1} \,
\end{eqnarray*}
and
\begin{eqnarray*}
\int_0^{\varepsilon} r^{\gamma_1}\int_{\RR^d}|\cF\Gamma(r)(\xi)|^2\mu(d\xi)dr= C \int_0^{\varepsilon} r^{\gamma_1 + \sum_{j=1}^d H_j-d} dr = C \varepsilon ^{\sum_{j=1}^d H_j-d+1+\gamma_1}\,.
\end{eqnarray*}
So  $({\bf H2})$ is satisfied with $\beta_1=\frac{\gamma_2}{2}+\sum_{j=1}^d H_j -d+1$ and $\beta_2=\sum_{j=1}^d H_j-d+1+\gamma_1$.
\rule{0.5em}{0.5em}
\\
\\

\textbf{Acknowledgement}: We would like to gratefully thank Dr. Jingyu Huang for useful discussion.

\end{document}